\RequirePackage{ifpdf}
\ifpdf % We are running pdfTeX in pdf mode
\documentclass[pdftex]{sigma}
\else
\documentclass{sigma}
\fi

\usepackage{amsmath}
\usepackage{amsfonts}
\usepackage{amssymb}
\usepackage[mathscr]{eucal} % Fuente Euscript, Euler Script font.
\usepackage{graphicx}
\usepackage{hyperref}
\usepackage[new]{old-arrows}
\usepackage{enumitem}
\usepackage{comment}
\usepackage{todonotes}
\usepackage[all]{xy}

\newcommand{\T}{\mathsf{T}}
\newcommand{\dd}{\mathrm{d}}

\usepackage[symbol]{footmisc}

\usepackage{amsthm}
\usepackage{thmtools}

\numberwithin{equation}{section}

\begin{document}

\renewcommand{\PaperNumber}{***}

\thispagestyle{empty}

%\FirstPageHeading

\ArticleName{On Gauge Poisson Brackets with Prescribed Symmetry}
\ShortArticleName{On Gauge Poisson Brackets with Prescibed Symmetry}

\Author{M. Avenda\~{n}o-Camacho\,${}^{\dag^{1}}$,
J. C. Ru\'iz-Pantale\'on\,${}^{\dag^{2}}$
and Yu. Vorobiev\,${}^{\dag^{3}}$}
\AuthorNameForHeading{Avenda\~{n}o-Camacho, Ru\'iz-Pantale\'on, Vorobiev}

\Address{$^{\dag}$~Departamento de Matem\'aticas, Universidad de Sonora, M\'exico}
\EmailD{$^{1}$\href{mailto:email@address}{misael.avendano\,@unison.mx},
        $^{2}$\href{mailto:email@address}{jose.ruiz\,@unison.mx},
        $^{3}$\href{mailto:email@address}{yury.vorobev\,@unison.mx}}
%\URLaddressD{\url{http://www.home.org/~mitome/}}

%\ArticleDates{Received ???, in final form ????; Published online ????}

\Abstract{In the context of averaging method, we describe a reconstruction of invariant connection-dependent Poisson structures from canonical actions of compact Lie groups on fibered phase spaces. Some symmetry properties of Wong's type equations are derived from the main results.}


%\Classification{17B63, 17B60, 53D17, 16W25}

    \section{Introduction}\label{sec:introdution}

A class of the so-called gauge Lie-Poisson brackets \cite{MoMaRa-84} naturally appears on an associated vector bundle of an $H$-principle bundle $H \rightarrow P \rightarrow Q$ and plays an important role in a Hamiltonian formulation of Wong's equations describing the motion of a ``colored'' particle in the Yang-Mills field \cite{Wong,Mon-84,Land-98,Maspfuhl-04,Holten-07,Sternberg-1977,Weinstein-1978}. The gauge Lie-Poisson brackets live on the pull-back $E=\tau^{\ast}(\mathrm{ad}^{\ast}(P))$ of the coadjoint bundle $\mathrm{ad}^{\ast}(P)=P\times_{H}\mathfrak{h}^{\ast}$ via the canonical projection $\tau:\T^{\ast}Q\rightarrow Q$ and are parameterized by $H$-invariant connections on the principle bundle $P$ \cite{GHV}. Then, on the phase space $E$ equipped with a gauge Lie-Poisson bracket, the Hamiltonian system with Hamiltonian $\mathcal{H}=\frac{1}{2}g^{ij}(q)p_{i}p_{j}$ (here $g=(g_{ij}(q))$ is a metric tensor on $Q$) just gives Wong's equations on $E$ \cite{Mon-84}.

The existence of first integrals of Wong's equations was studied, for example, in \cite{Holten-07, HorvNG-06}. Our point is to describe Wong's equations with prescribe symmetry, thinking of $\pi: E\rightarrow \T^{\ast}Q$ as a Lie-Poisson bundle over the symplectic base $(\T^{\ast}Q, \omega=\dd p_{i}\wedge \dd q^{i})$ and starting with an action of a compact connected Lie group $G$ on $E$ which leaves the fibers invariant and acts canonically with respect to the fiberwise Lie-Poisson structure \cite{Va-98}. This action may be not necessarily related with the $H$-action. Then, the question is reduced to the construction of $G$-invariant Lie-Poisson connections on $E$ and the reconstruction of the corresponding $G$-invariant Poisson brackets \cite{GHV,Va-98}. In this paper, we are interested in a more general setting of the above problem for a wide class of connection-dependent Hamiltonian systems on fibered spaces with symmetry which generalize Wong's equations. Our approach is based on the averaging method for Poisson connections \cite{MaMoRa-90,MHVV} and coupling Poisson structures on fibered spaces \cite{Vo-01,Vo-04,MYu-2012,MYu-17}. As a consequence of the main results, we get the following fact.

\paragraph{Claim 1.} Let $Q\times\mathfrak{h}^{\ast}$ be the product of a (pseudo)-Riemannian manidold $Q$ and the co-algebra $\mathfrak{h}^{\ast}$ equipped with the Lie-Poisson bracket \cite{Va-98}. Suppose that a compact connected  Lie group $G$ acts along the fibers of $Q\times\mathfrak{h}^{\ast}$ canonically with momentum map. Then, there exists a Yang-Mills gauge potential $\mathbf{A}\in\Omega^{1}(Q;\mathfrak{h})$, which induces a $G$-invariant Wong's equations on the phase space $\T^{\ast}Q\times\mathfrak{h}^{\ast}$. \\

As an application, we use this result to derive a gauge potential of the Wu-Yang potential \cite{Holten-07,HorvNG-06} by using the averaging with respect to an $\mathbb{S}^{1}$-action on the trivial $\mathfrak{so}^{\ast}(3)$-Poisson bundle.
%As application of the main result, we derive a gauge potential of the Wu-Yang potential \cite{Holten-07} by using the averaging with respect to an $\mathbb{S}^{1}$-action on the trivial $\mathrm{so}^{\ast}(3)$-Poisson bundle.

Moreover, we apply the main results to some class of generalized Wong equations on non-trivial Lie-Poisson bundles equipped with a class of fiberwise linear gauge Poisson brackets associated to transitive Lie algebroids \cite{Vo-04}.

\section{Nonlinear Gauge Poisson Brackets}

Let $\T^{\ast}Q$ be the cotangent bundle over a manifold $Q$ and $\omega=\dd p_{i}\wedge \dd q^{i}$ the canonical symplectic form. Suppose we are given a Poisson manifold $N$ equipped with a Poisson bracket
    \begin{equation*}
        \{y^{\alpha},y^{\beta}\}_{N}=\Psi^{\alpha\beta}(y),
    \end{equation*}
where $\Psi=\frac{1}{2}\Psi^{\alpha\beta}(y)\frac{\partial}{\partial y^{\alpha}}\wedge\frac{\partial}{\partial y^{\beta}}$ is the corresponding Poisson tensor \cite{Va-98}.

Assume we have an arbitrary 1-form on $Q\times N$ of the form
    \begin{equation*}
        \mathcal{A}=\mathcal{A}_{i}(q,y)\dd q^{i},
    \end{equation*}
which will be called a\textit{ gauge form}. Then one can associate to $\Psi$ and $\mathcal{A}$ the following 2-form $\mathcal{F}=\frac{1}{2}\mathcal{F}_{ij}(q,y)\dd q^{i}\wedge \dd q^{j}$ on $Q\times N$ with coefficients
    \begin{equation}
        \mathcal{F}_{ij}:=\frac{\partial\mathcal{A}_{j}}{\partial q^{i}}%
        -\frac{\partial\mathcal{A}_{i}}{\partial q^{j}}+\Psi^{\alpha\beta}%
        \frac{\partial\mathcal{A}_{i}}{\partial y^{\alpha}}\frac{\partial
        \mathcal{A}_{j}}{\partial y^{\beta}}, \label{Re1}%
    \end{equation}
or, in free coordinate notations,
    \begin{equation}
        \mathcal{F}=\dd_{Q}\mathcal{A+}\frac{1}{2}\{\mathcal{A}\wedge\mathcal{A}\}_{N}.
        \label{Re2}%
    \end{equation}
Here, we denote by $\dd_{Q}$ the partial exterior derivative along $Q$ for forms on $Q\times N$.
%Here, we denote by $d_{Q}$ and $d_{N}$ the partial exterior derivatives of forms on $Q\times N$ along $Q$ and $N$, respectively. So, $d=d_{Q}+d_{N}$ is the exterior derivative on $Q\times N$ and hence, we have $d_{Q}^{2}=d_{N}^{2}=0$ and $d_{Q}\circ d_{N}=-d_{N}\circ d_{Q}$.
%Denote by $\mathrm{Casim}_{Q}(N,\Psi)$ the space of smooth functions on $Q\times N$ which take values in the space of Casimir functions of $\Psi$ for each fixed $q$. \todo{Porqu\'e definir Casim aqu\'i?}

\begin{theorem}
Every gauge 1-form $\mathcal{A}$ induces a Poisson structure $\Pi_{\mathcal{A}}$ on the phase space $E=\T^{\ast}Q\times N$ given by the following pairwise bracket relations:
    \begin{align}
        &\{p_{i},p_{j}\} = \mathcal{F}_{ij}, \qquad\qquad\ \{p_{i},q^{j}\}=\delta_{i}^{j}, \quad \{q^{i},q^{j}\}=0, \label{YaM1} \\
        &\{p_{i},y^{\alpha}\} = -\Psi^{\alpha\beta}\frac{\partial\mathcal{A}_{i}}{\partial y^{\beta}}, \quad \{q^{i},y^{\alpha}\}=0, \label{YaM3} \\
        & \{y^{\alpha},y^{\beta}\} = \Psi^{\alpha\beta}. \label{YaM4}
    \end{align}
\end{theorem}

The verification of the Jacobi identity for these brackets is a direct computation by using formulas %relations
\eqref{Re1}-\eqref{Re2} and the Jacobi identity for $\Psi$.

We say that the bracket relations \eqref{YaM1}-\eqref{YaM4} define a \textit{gauge Poisson bracket} on $E$ associated to a gauge 1-form $\mathcal{A}$ or, shortly, an $\mathcal{A}$-\textit{Poisson bracket}.

The Poisson tensor $\Pi_{\mathcal{A}}$ corresponding to relations \eqref{YaM1}-\eqref{YaM4} is represented as
    \begin{equation}
        \Pi_{\mathcal{A}}:=\frac{1}{2}\Psi^{\alpha\beta}\frac{\partial}{\partial
        y^{\alpha}}\wedge\frac{\partial}{\partial y^{\beta}}+\frac{\partial}{\partial
        p_{i}}\wedge\mathrm{hor}_{i}+\frac{1}{2}\mathcal{F}_{ij}\frac{\partial
        }{\partial p_{i}}\wedge\frac{\partial}{\partial p_{j}}. \label{PT1}%
    \end{equation}
Here,
\begin{equation}
\mathrm{hor}_{i}=\mathrm{hor}_{i}^{\mathcal{A}}:=\frac{\partial
}{\partial q^{i}}+\Psi^{\alpha\beta}\frac{\partial\mathcal{A}_{i}}{\partial
y^{\alpha}}\frac{\partial}{\partial y^{\beta}} \label{Hor}%
\end{equation}
%denote the horizontal vector fields on $E$, \redcolor{which generates the horizontal distribution}
denote the vector fields generating the distribution $\mathbb{H}$ on $Q \times N$:
    \begin{equation}\label{Ho2}
        \mathbb{H}_{(q,y)}:=\mathrm{Span}\{ \mathrm{hor}_{i} \mid i=1,\ldots,\dim Q \}. % \subset T_{(q,y)}(Q \times N).}
    \end{equation}

Note also that at every point $(p,q,y)\in E$, we have
    \begin{equation*}
        \operatorname{rank}{\Pi}_{\mathcal{A}}=2\dim Q+\operatorname{rank}\Psi(y).
    \end{equation*}
Moreover, if $S\subset E$ is a symplectic leaf of $\Pi_{\mathcal{A}}$ through $(p,q,y)$, then
    \begin{equation}
        \T_{(p,q,y)}S = \overline{\mathcal{\mathbb{H}}}_{(p,q,y)}\oplus \T_{y}\mathcal{O}, \label{Ho1}%
    \end{equation}
where $\mathcal{O}$ is the symplectic leaf of $\Psi$ through $y\in N$ and
    \begin{equation}
        \mathbb{\overline{H}}_{(p,q,y)}:=\mathbb{H}_{(q,y)}\oplus\mathrm{Span}\left\{ \left. \frac{\partial}{\partial p_{i}} \,\right|\, i=1,\ldots,\dim Q \right\}. \label{EC1}
        %\redcolor{\tilde{\mathbb{H}}_{q,y}:=\mathrm{Span}\left\{ \left. \frac{\partial}{\partial p^{i}},\, \mathrm{hor}_{i} \,\right|\, i=1,\ldots,\dim Q \right\} \subset T_{(p,q,y)}E} % \label{Ho2}%
    \end{equation}
%\redcolor{is the horizontal distribution induced by $\mathbb{H}$ on $E$.}

%We say that bracket relations \eqref{YaM1}-\eqref{YaM4} define a \textit{gauge Poisson bracket} on $E=\T Q\times N$ associated to a gauge 1-form $\mathcal{A}$ or shortly, an $\mathcal{A}$-\textit{Poisson bracket}. \\

\paragraph{Geometric Interpretation.} Consider the product manifold $Q\times N$ as the total space of a trivial Poisson bundle over $Q$ with typical fiber $(N,\Psi)$ \cite{Va-98}. Then, each gauge form $\mathcal{A}$ determines an \textit{Ehresmann-Poisson connection} on $Q\times N$ associated to the horizontal distribution $\mathbb{H}$ in \eqref{Ho2} \cite{GHV,Va-98}. The curvature
    \begin{equation*}
        \mathcal{C}^{\mathbb{H}}=\frac{1}{2}C_{ij}^{\alpha}(q,y)\dd q^{i}\wedge
\dd q^{j}\otimes\frac{\partial}{\partial y^{\alpha}}%
    \end{equation*}
of this connection is defined by the 2-form $\mathcal{F}$ in \eqref{Re1},
    \begin{equation*}
        C_{ij}^{\alpha}=-\Psi^{\alpha\beta}\frac{\partial\mathcal{F}_{ij}}{\partial
y^{\beta}}.
    \end{equation*}
The vanishing of $\mathcal{C}^{\mathbb{H}}$ means the integrability of the horizontal distribution $\mathbb{H}$.

Consider the pullback $E=\T^{\ast}Q\times N=\tau^{\ast}(Q\times N)$ of the trivial Poisson bundle $Q\times N \rightarrow Q$ by the canonical projection $\tau:\T^{\ast}Q \rightarrow Q$. The cotangent bundle $\T^{\ast}Q$ is equipped with the standard symplectic form $\omega=\dd p_{i}\wedge \dd q^{i}$ and $\tau$ is a Lagrangian submersion. Therefore, $\pi_{1}:E\rightarrow \T^{\ast}Q$ is a trivial Poisson bundle over a symplectic manifold which carries an Ehresmann-Poisson connection $\mathbb{\overline{H}}$ \eqref{EC1} defined as the pull-back of $\mathbb{H}$ by $\tau$. Then, the $\mathcal{A}$-Poisson tensor $\Pi_{\mathcal{A}}$ in \eqref{PT1} is just the \textit{coupling Poisson structure} associated to the geometric data $(\Psi, \mathbb{\overline{H}}, \pi_{1}^{\ast}\omega-\overline{\tau}^{\ast}\mathcal{F})$, where $\overline{\tau}=\tau \times \mathrm{id}_{N}:E \rightarrow Q \times N$ \cite{Vo-01}. %\cite{Vo-01}.

Denote by $\mathrm{Casim}_{Q}(N,\Psi)$ the space of smooth functions $C=C(q,y)$ on $Q\times N$ which take values in the space of Casimir functions of $\Psi$ for each fixed $q$. If $C \in \mathrm{Casim}_{Q}(N,\Psi)$, then $\overline{\tau}^{\ast}C$ is a Casimir function of $\Pi_{\mathcal{A}}$ under the condition $\dd_{Q}C=0$.

\section{Hamiltonian Dynamics}

For every smooth function $\mathcal{H}=\mathcal{H}(p,q,y)$ on $\T^{\ast}Q\times N$, the Hamiltonian dynamics relative to a given $\mathcal{A}$-Poisson bracket and the Hamiltonian $\mathcal{H}$ is described by the system:
    \begin{align}
        \dot{p}_{j} &= -\frac{\partial\mathcal{H}}{\partial q^{j}}-\Psi^{\alpha\beta}\frac{\partial\mathcal{A}_{j}}{\partial y^{\alpha}}\frac{\partial\mathcal{H}}{\partial y^{\beta}}+\mathcal{F}_{ij}\frac{\partial\mathcal{H}}{\partial p_{i}}, \label{H1} \\
        \dot{q}^{j} &= \frac{\partial\mathcal{H}}{\partial p_{j}}, \label{H2} \\
        \dot{y}^{\beta} &=\Psi^{\alpha\beta}\frac{\partial\mathcal{H}}{\partial y^{\alpha}}+\frac{\partial\mathcal{H}}{\partial p_{i}}\Psi^{\alpha\beta}\frac{\partial\mathcal{A}_{i}}{\partial y^{\alpha}}. \label{H3}%
    \end{align}

In the case when $N$ is a co-algebra equipped with the Lie-Poisson structure \cite{Va-98}, we have the following class of Hamiltonian systems associated to Wong's equations.

Let $N=\mathfrak{h}^{\ast} \simeq \mathbb{R}_{y}^{n}$  be the co-algebra of a Lie algebra $\mathfrak{h}$ equipped with the Lie-Poisson bracket
    \begin{equation*}
        \{y^{\alpha},y^{\beta}\}_{\mathfrak{h}}=\lambda_{\gamma}^{\alpha\beta}y^{\gamma},
    \end{equation*}
where $\lambda_{\gamma}^{\alpha\beta}$ are the structure constants relative to a basis $\{e_{\alpha}\}$ of $\mathfrak{h}$. So, in this case, the Poisson tensor $\Psi$ has the components $\Psi^{\alpha\beta}=\lambda_{\gamma}^{\alpha\beta}y^{\gamma}$. Suppose we are given a \textit{vector gauge potential}, that is, an $\mathfrak{h}$-valued 1-form $\mathbf{A}\in\Omega^{1}(Q;\mathfrak{h})$ on $Q$,
    \begin{equation*}
        \mathbf{A}=A_{\alpha i}(q)\dd q^{i}\otimes e_{\alpha},
    \end{equation*}
and the corresponding \textit{Yang-Mills field}, defined as the $\mathfrak{h}$-valued 2-form on $Q$:
    \begin{equation}
        \mathbf{F:}=\dd \mathbf{A}+\frac{1}{2}[\mathbf{A}\wedge\mathbf{A}]_{\mathfrak{h}}=\frac{1}{2}F_{\alpha ij}(q)\dd q^{i}\wedge \dd q^{j}\otimes e_{\alpha}. \label{Cur}%
    \end{equation}
%In this case, a gauge form $\mathcal{A}=\mathcal{A}_{i}(q,y)dq^{i}$ and 2-form $\mathcal{F}=\frac{1}{2}\mathcal{F}_{ij}$ $(q,y)dq^{i} \wedge dq^{i}$ on $Q\times\mathfrak{h}^{\ast}$ are given by
%    \begin{equation*}
%        \mathcal{A}_{i}(q,y):=y^{\alpha}A_{\alpha i}(q), \quad \mathcal{F}_{ij}(q,y)=y^{\alpha}F_{\alpha ij}(q).
%    \end{equation*}
The corresponding gauge form $\mathcal{A}=\mathcal{A}_{i}(q,y)\dd q^{i}$ and 2-form $\mathcal{F}=\frac{1}{2}\mathcal{F}_{ij}(q,y)\dd q^{i}\wedge \dd q^{i}$ on $Q\times\mathfrak{h}^{\ast}$ are given by
    \begin{equation}
        \mathcal{A}_{i}(q,y):=y^{\alpha}A_{\alpha i}(q), \quad
        \mathcal{F}_{ij}(q,y)=y^{\alpha}F_{\alpha ij}(q). \label{Cur2}%
    \end{equation}
Therefore, in this case, the $\mathcal{A}$-Poisson bracket on the phase space $\T^{\ast}Q\times\mathfrak{h}^{\ast}$ is just a gauge Lie-Poisson bracket which is given by the relations
    \begin{align}
        &\{p_{i},p_{j}\} = y^{\alpha}F_{\alpha ij}(q), \qquad\quad\ \{p_{i},q^{j}\}=\delta_{i}^{j}, \quad \{q^{i},q^{j}\}=0, \label{LP1} \\
        &\{p_{i},y^{\alpha}\} = -\lambda_{\gamma}^{\alpha\beta}y^{\gamma}A_{\beta i}(q), \quad \{q^{i},y^{\alpha}\}=0, \label{LP2} \\
        &\{y^{\alpha},y^{\beta}\} = \lambda_{\gamma}^{\alpha\beta}y^{\gamma}. \label{LP3}
    \end{align}

Let $g=(g_{ij})$ be a pseudo-Riemannian metric on $Q$. The Hamiltonian dynamics \eqref{H1}-\eqref{H3} on $E=\T^{\ast}Q\times\mathfrak{h}^{\ast}$ relative to these brackets and the Hamiltonian
    \begin{equation}
        \mathcal{H}=\frac{1}{2}g^{ij}(q)p_{i}p_{j} \label{W1}%
    \end{equation}
takes the form of the classical Wong's equations,
    \begin{align}
        \dot{p}_{i} &= -\frac{1}{2}\frac{\partial g^{jk}}{\partial q^{i}}p_{j}p_{k}-y^{\alpha}F_{\alpha ij}g^{jk}(q)p_{k}, \label{W2} \\
        \dot{q}^{i} &= g^{ij}(q)p_{j}, \label{W3} \\
        \dot{y}^{\alpha} &=-y^{\gamma}\lambda_{\gamma}^{\alpha\beta}A_{\beta i}g^{ij}(q)p_{j}, \label{W4}
    \end{align}
determined by the vector gauge potential $\mathbf{A}$ and the corresponding Yang-Mills field $\mathbf{F}$.

Consider the Poisson tensor of \eqref{LP1}-\eqref{LP3} on $E=\T^{\ast}Q\times\mathfrak{h}^{\ast}$:
    \begin{equation}
    \Pi_{\mathbf{A}}:=\frac{1}{2}\lambda_{\gamma}^{\alpha\beta}y^{\gamma}%
    \frac{\partial}{\partial y^{\alpha}}\wedge\frac{\partial}{\partial y^{\beta}%
    }+\frac{\partial}{\partial p_{i}}\wedge\mathrm{hor}_{i}+\frac{1}%
    {2}y^{\alpha}F_{\alpha ij}(q)\frac{\partial}{\partial p_{i}}\wedge
    \frac{\partial}{\partial p_{j}}, \label{W5}%
    \end{equation}
where
    \begin{equation}
    \mathrm{hor}_{i}=\frac{\partial}{\partial q^{i}}+\lambda_{\gamma
    }^{\alpha\beta}y^{\gamma}A_{\alpha i}(q)\frac{\partial}{\partial y^{\beta}}
    \label{W6}%
    \end{equation}
can be interpreted as a set of \textit{linear vector fields} on the \textit{trivial} \textit{Lie-Poisson bundle} $Q\times\mathfrak{h}^{\ast}\rightarrow Q$ which generate the horizontal distribution associated to the Lie-Poisson connection $\nabla^{\mathbf{A}}$ given by
    \begin{equation*}
        \nabla^{\mathbf{A}}\xi = \dd\xi - \mathrm{ad}_{\mathbf{A}}^{\ast}\xi,
    \end{equation*}
for any smooth vector function $\xi:Q\rightarrow\mathfrak{h}^{\ast}$ determining a section of the trivial bundle \cite{GHV}. Therefore, under a fixed pseudo-metric $g$, the gauge Lie-Poisson structure $\Pi_{\mathbf{A}}$ and the classical Wong's equations are uniquely determined by $\nabla^{\mathbf{A}}$. The curvature of this connection is related with the $\mathfrak{h}$-valued 2-form in \eqref{Cur} by
    \begin{equation*}
        \mathrm{Curv}^{\nabla^{\mathbf{A}}}=-\mathrm{ad}_{\mathbf{F}}^{\ast}.
    \end{equation*}

%%%%%%%
Let $(p(t),q(t),y(t))$ be the trajectory of \eqref{W2}-\eqref{W4} passing through a point $(p^{0},q^{0},y^{0})$ at $t=0$. Then the last equation \eqref{W4} together with \eqref{W3} describes the $\nabla^{\mathbf{A}}$-parallel transport from $(q^{0},y^{0})$ along the curve $q(t)$, that is, $\nabla_{\dot{q}}^{\mathbf{A}}\dot{y}=0$ or, equivalently,
    \begin{equation}
        \dot{y}=\mathrm{ad}_{\mathbf{A}(\dot{q})}^{\ast}y. \label{CO}%
    \end{equation}
Let $H$ be the connected Lie group associated with $\mathfrak{h}$ and $\mathrm{Ad}_{\alpha}^{\ast}:\mathfrak{h}^{\ast}\rightarrow\mathfrak{h}^{\ast}$ the co-adjoint representation. Then, it follows from \eqref{CO} that the trajectory $(p(t),q(t),y(t))$ belongs to the product $\T^{\ast}Q\times\mathcal{O}$, where $\mathcal{O}\subset\mathfrak{h}^{\ast}$ is the co-adjoint orbit through $y^{0}$.

%\begin{Remark}
%%%%%%%%%%%%%%%%
%In the above derivation of Wong equations, we deal with the adjoint bundle $Q\times\mathfrak{h}^{\ast}$ of the trivial principle
%bundle $P=Q\times H$. In the general case, starting with a nontrivial principle bundle $H\rightarrow P \overset{\nu}{\rightarrow} Q$, first, we fix a Lie-Poisson connection $\nabla$ on the coadjoint bundle $\mathrm{ad}(P)=P\times_{H}\mathfrak{h}^{\ast}$ and then construct the $\nabla$-dependent gauge Lie-Poisson structure. Therefore, in local coordinates, the corresponding Wong's equations on the phase space $E=\tau^{\ast}\mathrm{ad}^{\ast}(P)$ can be written as \eqref{W2}-\eqref{W4}. A further generalization of this classical situation is related with a class of coupling Lie-Poisson structures coming from transitive Lie algebroids \cite{Vo-04}.
%%%%%%%%%%%%%%%
%\end{Remark}

\begin{example}
(Wu-Yang Monopole). Consider the case when $(Q=\mathbb{R}_{0}^{3}=\mathbb{R}^{3}\setminus\{0\},g=(g_{ij}(q)))$ and $\mathfrak{h}^{\ast}=\mathfrak{so}^{\ast}(3)$. On the phase space $\T^{\ast}\mathbb{R}_{0}^{3}\times\mathfrak{so}^{\ast}(3)$, consider a vector gauge potential $\mathbf{A}_{i}=\mathbf{A}_{\alpha i}(q) \otimes e_{\alpha} \in \mathbb{R}^{3}$ corresponding to the Wu-Yang monopole \cite{HorvNG-06, Holten-07},
    \begin{equation}
        \mathbf{A}_{i}(q)=\frac{1}{|q|^{2}}\,q\times e_{i} \quad (i=1,2,3), \label{WYP}%
    \end{equation}
where $\{e_\alpha\}$ is the canonical basis in $\mathbb{R}^{3}$. The corresponding Yang-Mills field $\mathbf{F}$ has the components
    \begin{equation*}
        F_{\alpha ij} = \frac{1}{|q|^{4}} \epsilon_{ijk}q^{k}q^{\alpha}.
        %= \operatorname{rot}{\mathbf{A}_{\alpha}} + \frac{1}{2}\epsilon^{\alpha\beta\gamma}\mathbf{A}_{\beta}\times\mathbf{A}_{\gamma} \\
        %= -\frac{q_{\alpha}}{|q|^{4}}\,q.
    \end{equation*}
In vector form, the corresponding Wong's type equations \eqref{H1}-\eqref{H3} on $\T^{\ast}\mathbb{R}_{0}^{3}\times\mathfrak{so}^{\ast}(3)$ are written as
    \begin{align}
        \frac{\dd q}{\dd t} &= \frac{\partial\mathcal{H}}{\partial p},\label{WY1} \\
        \frac{\dd p}{\dd t} &= -\frac{\partial\mathcal{H}}{\partial q}+\frac{y \boldsymbol{\cdot} q}{|q|^{4}}\,q\times\frac{\partial\mathcal{H}}{\partial p},\label{WY2} \\
        \frac{\dd y}{\dd t} & =\frac{1}{|q|^{2}}\left( q\times\frac{\partial\mathcal{H}}{\partial p} \right)\times y.\label{WY3}
    \end{align}
\end{example}

\begin{remark}
The connection forms on a trivial principle bundle $P=Q\times H$ are just parameterized by the gauge 1-forms $\mathbf{A}\in\Omega^{1}(Q;\mathfrak{g})$ \cite{GHV}. In general, if we start with a nontrivial principle bundle $H\rightarrow P$
$\rightarrow$ $Q$, then the connection forms on $P$ induce the linear connections $\nabla$ on the associated co-adjoint bundle $\mathrm{ad}^{\ast}(P)=P\times_{H}\mathfrak{h}^{\ast}$ which determine the gauge Lie-Poisson brackets on $\tau^{\ast}\mathrm{ad}^{\ast}(P)$ \cite{Mon-84, Land-98}. Locally, $\nabla=\nabla^{\mathbf{A}}$ and the corresponding Wong's equations on the phase space $E$ are written in the form \eqref{W2}-\eqref{W4}. A further generalization of this classical situation is related with a class of coupling Lie-Poisson structures coming from transitive Lie algebroids \cite{Vo-04} (see also Section \ref{sec:NonTrivialBundle}).
\end{remark}

\section{$G$-Symmetry}

Consider the product manifold $E_{Q}=Q\times N$ as a total space of the trivial Poisson bundle over $Q$ with typical fiber $(N,\Psi)$ \cite{Va-98}. Recall that $\mathrm{Casim}_{Q}(N,\Psi)\subset C^{\infty}(E_Q)$
 denotes the space of Casimir functions of the vertical Poisson structure on $E_{Q}$ induced by $\Psi$.
%which will be denoted in the same way.
Let $G$ be a compact connected Lie group and $\mathfrak{g}$ its Lie algebra. Denote by $\dd g$ the normalized Haar measure on $G$.

Assume that we are given a  smooth family of Hamiltonian $G$-actions on $E_{Q}$, that is, a smooth (left) $G$-action $\Phi:G\times E_{Q}\rightarrow E_{Q}$ on $E_Q$ such that each fiber is invariant under the diffeomorphism $\Phi_{g}:E_{Q}\rightarrow E_{Q}$,
    \begin{equation*}
        \Phi_{g}(q,y) = \big( q,\Phi_{g}^{q}(y) \big) \quad \forall g\in G,%
    \end{equation*}
where $(q,y)\in E_{Q}$ and $\Phi_{g}^{q}:N\rightarrow N$ is a diffeomorphism parameterized by $q\in Q$. Moreover, we assume that there exists a fiberwise momentum map $\mathbb{J}:E_{Q}\rightarrow\mathfrak{g}^{\ast}$, that is, for every $a\in\mathfrak{g}$, the infinitesimal generator is
    \begin{equation}
        \Upsilon_{a}:=\left. \frac{\dd}{\dd t}\big[ \Phi_{\exp(ta)}(q,y) \big]\right|_{t=0}=-\Psi^{\alpha\beta}\frac{\partial\mathbb{J}_{a}}{\partial y^{\beta}}(q,y)\frac{\partial}{\partial y^{\alpha}}. \label{IG1}%
    \end{equation}
Here, $\mathbb{J}_{a}(q,y)=\langle \mathbb{J}(q,y),a \rangle$. Therefore, the infinitesimal generator $\Upsilon_{a}$ of $\Phi$ coincides with the vertical Hamiltonian vector field of the function $\mathbb{J}_{a}:E_{Q}\rightarrow\mathbb{R}$. Since $E=\tau^{\ast}E_{Q}$, the $G$-action $\Phi$ is canonically lifted to $E$. This lifted action is also acting along the fibers of $E$ in a Hamiltonian fashion, where the fiberwise momentum map is the pull back of $\mathbb{J}$ by $\tau:\T^{\ast}Q\rightarrow Q$.

Now, let us define a gauge 1-form $\mathcal{A}=\mathcal{A}_{i}(q,y)\dd q^{i}$ on $E_{Q}$ by the formula
    \begin{equation}
        \mathcal{A}_{i}(q,y):=-\int\limits_{G}\left(
        \int\limits_{0}^{1} \frac{\partial\mathbb{J}_{a}}{\partial q^{i}} \big( q,\Phi_{\exp(ta)}^{q}(y) \big)\dd t \right) \dd g + C_{i}, \label{IGF1}%
    \end{equation}
where $g=\exp a$ for a certain $a\in\mathfrak{g}$ and $C_{i}\in\mathrm{Casim}_{Q}(N,\Psi)$ are arbitrary Casimir functions.
%Moreover,
%the integration over $G$ is taken with respect to the normalized Haar mesuare
%$dg$.
For a more detailed definition of the right hand side of \eqref{IGF1} see, for example, \cite{MaMoRa-90,JYu,MHVV}, where formula \eqref{IGF1} appears in the context of the so-called Hannay-Berry connection.

\begin{lemma}\label{lemhor}
The horizontal distribution $\mathbb{H}$ \eqref{Ho2} on $E_Q$ associated to a gauge form $\mathcal{A}$ is invariant with respect to the $G$-action if and only if $\mathcal{A}$ is defined by \eqref{IGF1} for some $C_{i}\in\mathrm{Casim}_{Q}(N,\Psi)$.
\end{lemma}

%%%%%%%%
The proof of this Lemma can be found in \cite{MHVV} which is based on the following arguments. Denote by $\langle \,\cdot\, \rangle_{G}=\int_{G}\Phi_{g}^{\ast}(\cdot)\dd g$ the $G$-average which is well-defined for tensor fields on $E_Q$ and $E$. Then, every Casimir function in $\mathrm{Casim}_{Q}(N,\Psi)$ is $G$-invariant. Moreover, the invariance of the horizontal distribution $\mathbb{H}$ with respect to the $G$-action is equivalent to the $G$-invariance of the generators $\mathrm{hor}_{i}$ in \eqref{Hor}, that is, $\langle\mathrm{hor}_{i}\rangle_{G}=\mathrm{hor}_{i}$, for $i=1,\ldots,\dim Q$. Equivalently, this condition, in terms of the infinitesimal generators $\Upsilon_{a}$ \eqref{IG1}, means that $\mathrm{hor}_{i}$ and $\Upsilon_{a}$ commute and can be rewritten as follows:
    \begin{equation}
        \frac{\partial\mathbb{J}_{a}}{\partial q^{i}}+\{\mathcal{A}_{i}%
        ,\mathbb{J}_{a}\}_{N}\in \mathrm{Casim}_{Q}(N,\Psi), \label{IC1}%
    \end{equation}
for all $a\in\mathfrak{g}$ and $i=1,\ldots,\dim Q$.
%The resolving of these
%relations relative to a gauge 1-form $\mathcal{A}$ leads to \eqref{IGF1}.
By the averaging arguments, \eqref{IC1} is equivalent to the following equations
 \begin{equation*}
     L_{\Upsilon_{a}}\mathcal{A}_i =- \frac{\partial\mathbb{J}_{a}}{\partial q^{i}} + \left\langle \frac{\partial\mathbb{J}_{a}}{\partial q^{i}} \right\rangle_{G},
 \end{equation*}
whose resolving leads to \eqref{IGF1} \cite{MHVV, MaMoRa-90, JYu}.

%%%%%%%

As a consequence of Lemma \ref{lemhor}, we derive also that the 2-form $\mathcal{F}$ in \eqref{Re2} associated to $\mathcal{A}$ in \eqref{IGF1} is also $G$-invariant,
    \begin{equation}
        \langle\mathcal{F}\rangle_{G}=\mathcal{F}\text{.}\label{IC}%
    \end{equation}
Indeed, averaging the both sides of the equality%
    \begin{equation*}
        \lbrack\mathrm{hor}_{i},\mathrm{hor}_{j}]=\Psi^{\alpha\beta}\frac{\partial\mathcal{F}_{ij}}{\partial y^{\alpha}}\frac{\partial}{\partial y^{\beta}}
    \end{equation*}
and taking into account the $G$-invariance of $\mathrm{hor}_{i}$ and $\Psi$, we get that $\langle\mathcal{F}_{ij}\rangle_{G}=\mathcal{F}_{ij}+K_{ij\text{ }}$ for some $K_{ij}\in\mathrm{Casim}_{Q}(N,\Psi)$. Consequently, $\langle K_{ij}\rangle_{G}=0$. So, from the $G$-invariance of the Casimir functions, we conclude that $K_{ij}=\langle K_{ij}\rangle_{G}=0$.

Recall that the $G$-action $\Phi_{q}$ on $E_{Q}$ is canonically lifted to a $G$-action $\overline{\Phi}_{q}$ on the pull-back bundle $E=\T^{\ast}Q\times N=\tau^{\ast}E_{Q}$, and that we denote $\overline{\tau}: E\rightarrow E_{Q}$, $\overline{\tau}=\tau\times\mathrm{id}_{N}$. We arrive at the main result.

\begin{theorem}\label{maintheo}
The gauge Poisson tensor $\Pi_{\mathcal{A}}$ associated to the 1-form $\mathcal{A}$ given by \eqref{IGF1} is $G$-invariant,%
    \begin{equation*}
        \overline{\Phi}_{g}^{\ast}\Pi_{\mathcal{A}}=\Pi_{\mathcal{A}} \quad \forall g\in G.
    \end{equation*}
The Hamiltonian system \eqref{H1}-\eqref{H3} with invariant Hamiltonian $\mathcal{H}$ is also $G$-invariant (in particular, with $\mathcal{H}$ of the form \eqref{W1}). Moreover, if the momentum map $\mathbb{J}_a$  satisfies the condition
    \begin{equation}
        \left\langle \frac{\partial\mathbb{J}_a}{\partial q_{j}} \right\rangle_{G}=0 \quad \forall a\in \mathfrak{g} \ \text{and} \ j, \label{AC}%
    \end{equation}
then
    \begin{enumerate}[label=(\alph*)]
      \item the lifted $G$-action on $E$ is canonical with respect to $\Pi_{\mathcal{A}}$ with momentum map $\overline{\tau}^{\ast}\mathbb{J}=E\rightarrow\mathfrak{g}^{\ast}$;

      \item the components of $\overline{\tau}^{\ast}\mathbb{J}$ are first integrals of motion of the Hamiltonian system \eqref{H1}-\eqref{H3},
          \begin{equation}
            \{\mathcal{H},\mathbb{J}_{a}\circ\overline{\tau}\}_{\mathcal{A}}=0 \quad \forall a\in\mathfrak{g.} \label{FI}%
          \end{equation}
    \end{enumerate}
\end{theorem}

The first part of the theorem follows from Lemma \ref{lemhor} and property \eqref{IC}. To verify property \eqref{FI}, we remark that
    \begin{equation*}
        L_{\mathrm{hor}_{j}^{\mathcal{A}}}\mathbb{J}_a=\left\langle \frac{\partial\mathbb{J}_a}{\partial q_{j}} \right\rangle_{G}%
    \end{equation*}
and hence by \eqref{AC}, we have $L_{\mathrm{hor}_{j}^{\mathcal{A}}}\mathbb{J}_a=0$. This, together with representations \eqref{PT1} and \eqref{Hor}, implies that the lifted $G$-action $\overline{\Phi}$ is also Hamiltonian with respect to $\Pi_{\mathcal{A}}$ with momentum map $\mathbb{J}\circ\overline{\tau}$ and the infinitesimal generators $\overline{\Upsilon}_{a}=\mathbf{i}_{\dd(\mathbb{J}_{a}\circ\tau)}\Pi_{\mathcal{A}}$. Finally, taking into account the $G$-invariance of $\mathcal{H}$, we get $0=L_{\overline{\Upsilon}_{a}}\mathcal{H}=\{\mathbb{J}_{a}\circ\overline{\tau},\mathcal{H}\}_{\mathcal{A}}$.

\begin{remark}
Since the Lie group is compact, one can choose Casimir functions $C_{i}$ in \eqref{IGF1} in such way that the momentum map $\mathbb{J}$ becomes equivariant (see, \cite{OrtRa-04}).
\end{remark}

\begin{corollary}
In the case when $N=\mathfrak{h}^{\ast}$ is a coalgebra equipped with the Lie-Poisson structure, the gauge 1-form $\mathcal{A}$ given by formula \eqref{IGF1} induces Wong's equations \eqref{W2}-\eqref{W4} on $\T^{\ast}Q\times\mathfrak{h}^{\ast}$ with prescribed $G$-symmetry. Under condition \eqref{AC}, the components of the momentum map  are the first integrals of Wong's equations.
\end{corollary}
Observe that Claim 1 in the introduction follows from this corollary.

\section{Particular Cases}

\paragraph{The case of $\mathbb{S}^{1}$-action.} Suppose we are given a vertical Hamiltonian vector field $\Upsilon=-\Psi^{\alpha\beta}\frac{\partial\mathbb{J}}{\partial y^{\beta}}(q,y)\frac{\partial}{\partial y^{\alpha}}$ of a smooth function $\mathbb{J}$
%\in C^{\infty}(E_{Q})
on $E_Q=Q\times N$ whose flow $\mathrm{Fl}_{\Upsilon}^{t}$ is $2\pi$-periodic, that is, $\mathrm{Fl}_{\Upsilon}^{t+2\pi}(m)=\mathrm{Fl}_{\Upsilon}^{t}(m)$ for any $t\in\mathbb{R}$ and $m=(q,y)\in E_{Q}$. This gives a family of Hamiltonian actions on $E_{Q}$ of the circle $\mathbb{S}^{1}=\mathbb{R}/2\pi\mathbb{Z}$ with Haar measure
$\frac{\dd t}{2\pi}$. Direct computation shows that the gauge form $\mathcal{A}=\mathcal{A}_{i}(q,y)\dd q^{i}$ in \eqref{IGF1} has the components: %is represented as follows
    \begin{equation*}
        \mathcal{A}_{i}(q,y)=\frac{1}{2\pi}\int\limits_{0}^{2\pi}(t-2\pi)\frac{\partial\mathbb{J}}{\partial q^{i}}\big( q,\mathrm{Fl}%
_{\Upsilon}^{t}(y) \big)\dd t + C_{i}.
    \end{equation*}
The normalization condition $\langle \mathcal{A}_{i} \rangle_{\mathbb{S}^{1}}=0$ holds if
    \begin{equation*}
        C_{i}=\pi\left\langle \frac{\partial \mathbb{J}}{\partial q^{i}} \right\rangle_{\mathbb{S}^{1}}=\frac{1}{2}\int\limits_{0}^{2\pi}\frac{\partial\mathbb{J}}{\partial q^{i}}\big( q,\mathrm{Fl}_{\Upsilon}^{t}(y) \big)\dd t.
    \end{equation*}
In this case, we get%
    \begin{equation}
        \mathcal{A}_{i}(q,y)=\frac{1}{2\pi}\int\limits_{0}^{2\pi}(t-\pi)\frac{\partial\mathbb{J}}{\partial q^{i}}\big( q,\mathrm{Fl}_{\Upsilon}^{t}(y) \big)\dd t. \label{PC}%
    \end{equation}

Consider again the special case when $N=\mathfrak{h}^{\ast}$ is a coalgebra equipped with the Lie-Poisson bracket \cite{Va-98}. Let $H$ be the connected Lie group associated with the Lie algebra $\mathfrak{h}$. Consider the co-adjoint action $\mathrm{Ad}_{h}^{\ast}:\mathfrak{h}^{\ast}\rightarrow\mathfrak{h}^{\ast}$ and assume that there exists a smooth section $s$ of the trivial bundle $Q\times\mathfrak{h}$, that is, an $\mathfrak{h}$-valued $C^{\infty}$-function $Q\ni q\mapsto s(q)\in\mathfrak{h}$ which satisfies the ``periodicity'' condition:
    \begin{equation}
        \exp(2\pi\,\mathrm{ad}_{s(q)})=\mathrm{id} \quad \forall q\in Q. \label{PC1}%
    \end{equation}
Then, we have the $\mathbb{S}^{1}$-action
    \begin{equation}
    (q,y)\mapsto(q,\mathrm{Ad}_{\exp ts(q)}^{\ast}y). \label{PC2}%
    \end{equation}
Taking into account that $\mathrm{Ad}_{\exp a}^{\ast}=\exp(\mathrm{ad}_{a}^{\ast})$, we get that the infinitesimal generator is a linear Hamiltonian vector field on $\mathfrak{h}^{\ast}$:
    \begin{equation}
        \Upsilon(q,y)=\mathrm{ad}_{s(q)}^{\ast}y\cdot\frac{\partial}{\partial
        y}\label{PC3}%
    \end{equation}
and hence
    \begin{equation}
        \mathbb{J}(q,y)=\langle s(q),y \rangle.\label{PC4}%
    \end{equation}
Formula \eqref{PC} for the components of the corresponding gauge 1-form $\mathcal{A} = \mathcal{A}_{i}(q,y)\dd q^{i}$ is written as follows
    \begin{equation}
        \mathcal{A}_{i}(q,y)=\frac{1}{2\pi}\int\limits_{0}^{2\pi}(t-\pi)
        \left\langle \exp(t\,\mathrm{ad}_{s(q)}^{\ast})y,\frac{\partial s}{\partial q^{i}}(q) \right\rangle \dd t. \label{PC5}%
    \end{equation}

Consider the case $\mathfrak{h}=\mathfrak{so}(3) \simeq \mathbb{R}^{3}$. The coalgebra $N=\mathfrak{so}^{\ast}(3) \simeq \mathbb{R}^{3}$ is equipped with the cyclic Lie-Poisson bracket
    \begin{equation*}
        \{y^{1},y^{2}\} = y^{3}, \quad \{y^{2},y^{3}\} = y^{1}, \quad \{y^{3},y^{1}\} = y^{2}.
    \end{equation*}
Fix a smooth mapping $s:Q \to \mathbb{R}^{3}$ with $|s(q)|=1$, $\forall q \in Q$. Then, condition \eqref{PC1} holds and the fiberwise Hamiltonian $\mathbb{S}^{1}$-action \eqref{PC2} on $E_{Q}=Q \times \mathfrak{so}^{\ast}(3)$ is just given by the rotations in the $\mathbb{R}^{3}_{y}$-space around the $s(q)$-axis. The corresponding infinitesimal generator and the momentum map take the form
    \begin{equation*}
        \Upsilon = (y \times s(q)) \cdot \frac{\partial}{\partial y}, \quad \mathbb{J}=\langle s(q),y \rangle.
    \end{equation*}
Then, by direct computation, we show that, in this case, formula \eqref{PC5} gives the following gauge 1-form
    \begin{equation}\label{eq:Aiqy}
        \mathcal{A}_{i}(q,y) = \left\langle s(q) \times y, \frac{\partial s}{\partial q^{i}}(q) \right\rangle.
    \end{equation}

Therefore, one can associated to \eqref{eq:Aiqy} the family of $\mathbb{S}^{1}$-invariant gauge Poisson brackets and Wong's equations on the phase space $\T^{\ast}Q \times \mathfrak{so}^{\ast}(3)$ which are parametrized by the smooth section $s:Q \to \mathbb{S}^{2}$.

\begin{example}
(The Wu-Yang monopole from an $\mathbb{S}^{1}$-action) Let $Q=\mathbb{R}_{0}^{3}$ and
    \begin{equation*}
        s(q) = \frac{q}{|q|}.
    \end{equation*}
Then, formula \eqref{eq:Aiqy} reads
    \begin{equation}\label{eq:Aqxy}
        \mathcal{A} = \frac{1}{|q|^{2}}(q \times y) \cdot\dd q
    \end{equation}
and gives the vector potential of the Wu-Yang monopole in \eqref{WYP}. As a consequence of Theorem \ref{maintheo}, the momentum map $\mathbb{J}=\langle \frac{q}{|q|},y \rangle$ is a first integral of Wong's equations \eqref{WY1}-\eqref{WY3} (see, also \cite{MoMaRa-84, Wong}). Finally, the $\mathbb{S}^{1}$-invariant gauge Poisson bracket associated to \eqref{eq:Aqxy} is given by the relations:
    \begin{align*}
         &\{p_{i},p_{j}\} = -\frac{1}{|q|^{4}}\epsilon_{ijk}q^{k}q^{\alpha}y^{\alpha}, \qquad \{p_{i},q^{j}\}= \delta_{i}^{j}, \quad \{q^{i},q^{j}\}=0, \\
         &\{p_{i},y^{\alpha}\} = -\frac{1}{|q|^{2}}\epsilon_{\alpha\beta\gamma}\epsilon_{\beta ij}q^{j}y^{\gamma}, \quad \{q^{i},y^{\alpha}\}=0, \\
         &\{y^{\alpha},y^{\beta}\} = -\epsilon_{\alpha\beta\gamma}y^{\gamma}.
    \end{align*}
\end{example}

\paragraph{The case of $\mathbb{T}^{r}$-action.} Suppose we have an action of the torus $G=\mathbb{T}^{r}=\mathbb{R}^{r}/
2\pi\mathbb{Z}^{r}$ with infinitesimal generators
    \begin{equation*}
        \Upsilon_{a}=a^{j}\Upsilon_{j},
    \end{equation*}
where $\{e_{j}\}$ is a basis of $\mathbb{R}^{r}$, $a=a^{j}e_{j}$, such that the flow $\mathrm{Fl}_{\Upsilon_{j}}^{t}$ of each $\Upsilon_{j}=\Upsilon_{e_{j}}$ is $2\pi$-periodic. Here, $\mathbb{J}_{a}(q,y)=a^{j}\mathbb{J}_{j}(q,y)$ with $\mathbb{J}_{j}=\mathbb{J}_{e_{j}}$. Then, formula \eqref{PC} reads
    \begin{equation*}
        \mathcal{A}_{i}(q,y):=-\int\limits_{0}^{2\pi}\cdots%
        \int\limits_{0}^{2\pi}
        \left(\int\limits_{0}^{1}a^{j}\frac{\partial\mathbb{J}_{j}}{\partial q^{i}} \big( q,\mathrm{Fl}%
        _{\Upsilon_{1}}^{a_{1}t}\circ\cdots\circ\mathrm{Fl}_{\Upsilon_{r}}^{a_{r}%
        t}(y) \big)\dd t \right) \dd a^{1} \cdots \dd a^{r}+C_{i}.%
    \end{equation*}

\paragraph{The case of $\mathrm{SO}(3)$-action.} Consider the basis $\{e_1,e_2,e_3 \}  $ of $\mathfrak{so}(3)$  satisfying the bracket relations
    \begin{equation*}
        [e_i,e_j] = \epsilon_{ijk} e_k.
    \end{equation*}
%given by
 %   \begin{equation*}
        %E_1=\begin{pmatrix}
        %0 & 0 & 0 \\
    %    0 & 0 &-1 \\
     %   0 & 1 & 0
      %  \end{pmatrix}, %\quad
        %E_2=\begin{pmatrix}
    %    0 & 0 & 1 \\
    %    0 & 0 & 0 \\
    %    -1 & 0 & 0
    %    \end{pmatrix}, \quad
        %E_3=\begin{pmatrix}
    %    0 & -1 & 0 \\
     %   1 & 0 & 0 \\
     %   0 & 0 & 0
     %   \end{pmatrix}.
%    \end{equation*}
Let $D^3_\pi =\{\mathbf{a}= a^ie_i \mid (a^1)^2+(a^2)^2+(a^3)^2 < \pi^2 \} \subset \mathfrak{so}(3)$. The exponential map $\exp: D^3_\pi\rightarrow \mathrm{SO}(3)$ is a diffeomorphism from $D^3_\pi$ onto its image and continuous over the closure of $D^3_\pi$. Assume we have an $\mathrm{SO}(3)$-Hamiltonian action on $E_Q$ with fiberwise momentum map $\mathbb{J}: E_Q \rightarrow \mathfrak{so}(3)$. For $\mathbf{a}= a^ie_i$ and $\mathbb{J}_\mathbf{a}(q,y)=a^j\mathbb{J}_j(q,y)$, with $\mathbb{J}_j=\mathbb{J}_{e_j} $. In this case, formula
\eqref{PC} for the gauge form is given by
    \begin{equation*}
        \mathcal{A}_{i}(q,y):=- \int_{| \mathbf{a}  | <\pi } \frac{\sin^2({|\mathbf{a}|}/{2})}{2 |\mathbf{a}|^2}\left( a^j \int_{0}^{1}\frac{\partial \mathbb{J}_j}{\partial q^i}\big(q, \Phi^q_{\exp t\mathbf{a}}(y) \big) \dd{t}\right)\dd a^1 \dd a^2 \dd a^3 +C_{i}.
    \end{equation*}

%%%%%%%%%%%%%%%%%%%%%%%%%%%%%%%%%%%%%%%%%%
%\section{Patents}

%This section is not mandatory, but may be added if there are patents resulting from the work reported in this manuscript.

%%%%%%%%%%%%%%%%%%%%%%%%%%%%%%%%%%%%%%%%%%
%\vspace{6pt}

%%%%%%%%%%%%%%%%%%%%%%%%%%%%%%%%%%%%%%%%%%
%% optional
%\supplementary{The following supporting information can be downloaded at:  \linksupplementary{s1}, Figure S1: title; Table S1: title; Video S1: title.}


\section{Nontrivial Lie-Poisson Bundles}\label{sec:NonTrivialBundle}

Let $(\pi_{Q}:E_{Q}\rightarrow Q,\Psi)$ be a Lie-Poisson bundle, that is, a vector bundle over $Q$ equipped with a vertical Poisson tensor $\Psi\in \Gamma(\wedge^{2}\T E_{Q})$ whose restriction to each fiber $E_{q}=\pi^{-1}(q)$ defines a linear Poisson structure (a fiberwise Lie-Poison structure) \cite{Va-98}:
    \begin{equation}
        \Psi_{q}=\lambda_{\gamma}^{\alpha\beta}(q)y^{\gamma}\frac{\partial}{\partial
        y^{\alpha}}\wedge\frac{\partial}{\partial y^{\beta}}.\label{VPSH}%
    \end{equation}
The dual bundle $E_{Q}^{\ast}$ becomes a bundle of Lie algebras with fiberwise Lie algebra structure $[\cdot,\cdot]_{\mathrm{fib}}:\Gamma(E_{Q}^{\ast})\times\Gamma(E_{Q}^{\ast})\rightarrow\Gamma(E_{Q}^{\ast})$ defined by%
    \begin{equation*}
        \ell_{\lbrack\eta_{1},\eta_{2}]_{\mathrm{fib}}}=\Psi(\dd\ell_{\eta_{1}},\dd\ell_{\eta_{2}}).
    \end{equation*}
Here, $\ell:\,\Gamma(E_{Q}^{\ast})\rightarrow C_{\mathrm{lin}}^{\infty}(E_{Q})$ is the natural identification between the space of smooth sections of $E_{Q}^{\ast}$ and the space of smooth fiberwise linear functions on $E_{Q}$. Denote by $\mathrm{ad}_{\eta}:\Gamma(E_{Q}^{\ast})\rightarrow\Gamma(E_{Q}^{\ast})$ the corresponding adjoint operator,
$\mathrm{ad}_{\eta}(\cdot)=[\eta,\cdot]_{\mathrm{fib}}$. Therefore, we also have the co-adjoint operator $\mathrm{ad}_{\eta}^{\ast}:\Gamma(E_{Q})\rightarrow\Gamma(E_{Q})$.

Assume that we are given a pair $(\nabla,\mathbf{F})$ consisting of a linear connection $\nabla$ on $E_{Q}$ \cite{GHV}
    \begin{equation*}
        \nabla:\mathfrak{X}(Q)\times\Gamma(E_{Q})\rightarrow\Gamma(E_{Q}),
    \end{equation*}
and a vector-valued 2-form on the base $\mathbf{F}\in\Omega^{2}(Q)\otimes\Gamma(E_{Q})$ satisfying the following conditions:
    \begin{enumerate}[label=(\roman*)]
        \item\label{CC1} $\nabla$ preserves the fiberwise Lie-Poisson bracket,
            \begin{equation}
                \lbrack\nabla_{v},\mathrm{ad}_{\eta}^{\ast}]=-\mathrm{ad}%
                _{\nabla_{v}^{\ast}\eta}^{\ast}, \label{LPVH1}%
            \end{equation}
        for any $v\in\mathfrak{X}(Q)$ and $\eta\in\Gamma(E_{Q}^{\ast})$;

        \item\label{CC2} the curvature 2-form $\mathrm{Curv}^{\nabla}\in\Omega^{2}(Q)\otimes\Gamma(\mathrm{End}(E_{Q}))$ takes values in the co-adjoint
        bundle,
            \begin{equation}
                \mathrm{Curv}^{\nabla}=-\mathrm{ad}^{\ast}_{\mathbf{F}};\label{LPVH2}%
            \end{equation}

        \item\label{CC3} $\mathbf{F}$ is $\nabla$-covariantly constant,
            \begin{equation}
                \nabla\mathbf{F}=0.\label{LPVH3}%
            \end{equation}
    \end{enumerate}

To every $v\in\mathfrak{X}(Q)$, we associate a linear vector field $\mathrm{hor}_{v}^{\nabla}$ on $E_{Q}$ uniquely defined by the condition
$L_{\mathrm{hor}_{v}^{\nabla}}\ell_{\eta}=\ell_{\nabla_{v}^{\ast}\eta}$, for any $\eta\in\Gamma(E_{Q}^{\ast})$. These vector fields generate the horizontal subbundle $\mathbb{H}^{\nabla}\subset \T E_{Q}$ so that we have the splitting
    \begin{equation}
        \T E_{Q}=\mathbb{H}^{\nabla}\oplus\mathbb{V},\label{SL}%
    \end{equation}
where $\mathbb{V}=\ker \dd\pi_{Q}$ is the vertical subbundle. Condition (\ref{LPVH1}) says that $\mathrm{hor}_{v}^{\nabla}$ is a
Poisson vector field of $\Psi$, $L_{\mathrm{hor}_{v}^{\nabla}}\Psi=0$.

Consider the pull-back bundle $E=\tau^{\ast}E_{Q}$ of $E_{Q}$ via the canonical projection $\tau:\T^{\ast}Q\rightarrow Q$. Therefore, $E$ is a Lie-Poisson bundle over the symplectic base $\T^{\ast}Q$. We have also the splitting $\T E = \overline{\mathbb{H}}^{\nabla} \oplus \overline{\mathbb{V}}$, which is the pullback of splitting (\ref{SL}) via $\tau$. According to the coupling method for Poisson structures, we have the following result \cite{Vo-01,Vo-04}.

\begin{proposition}
Every pair $(\nabla,\mathbf{F})$ satisfying conditions (\ref{LPVH1})-(\ref{LPVH3}) induces a Poisson tensor $\Pi^{\nabla,\mathbf{F}}$ on the total space $E$, which has the decomposition
    \begin{equation*}
        \Pi^{\nabla,\mathbf{F}}=\Pi_{H}^{\nabla,\mathbf{F}}+\Pi_{V}%
    \end{equation*}
into the horizontal and vertical components relative to the splitting (\ref{SL}). The bivector field $\Pi_{H}^{\nabla,\mathbf{F}}$ is tangent to the horizontal distribution $\overline{\mathbb{H}}^{\nabla}$ and defined by $\nabla$ and $\mathbf{F}$, the vertical part $\Pi_{V}$ is the canonical pull-back of the Poisson bivector field $\Psi$.
\end{proposition}

Remark that in coordinates $(p_{i},q^{j},y^{\alpha})$ on the total space of the Lie-Poisson bundle $E \to \T^{\ast}Q$, the pairwise Poisson brackets associated to $\Pi^{\nabla,\mathbf{F}}$ are just given by relations \eqref{LP1}-\eqref{LP3}, where the (local) gauge 1-form $\mathbf{A}$ on $Q$ with values in $\Gamma_{\mathrm{loc}}(E_{Q}^{\ast})$ is defined by $\nabla\zeta = \dd\zeta - \mathrm{ad}_{\mathbf{A}}^{\ast}\zeta$ for $\zeta \in \Gamma_{\mathrm{loc}}(E_{Q})$.

As consequence of proposition above, one can define (generalized) Wong's equations associated to the data $(\nabla,\mathbf{F})$ and a fixed pseudo-metric on $Q$ as the Hamiltonian system
    \begin{equation}
        \Big( E,\, \Pi^{\nabla,\mathbf{F}},\, \mathcal{H}=\frac{1}{2}g^{ij}(q)p_{i}p_{j} \Big). \label{NW}%
    \end{equation}

In particular, it is easy to see that, in the trivial case when $E_{Q}=Q\times\mathfrak{h}^{\ast}$ and $\nabla=\nabla^{\mathbf{A}}$ for some gauge potential $\mathbf{A}\in\Omega^{1}(Q;\mathfrak{h})$ and $\mathbf{F}$ given by \eqref{Cur}, conditions (\ref{LPVH1})-(\ref{LPVH3}) hold. This shows that $\Pi^{\nabla,\mathbf{F}}=\Pi_{\mathbf{A}}$.

As we mentioned above, a realization of the described scheme comes from transitive Lie algebroids.

\begin{example}
Let $(A,\rho,\{\,,\,\}_{A})$ be a transitive Lie algebroid over a manifold $Q$, where $\rho:A\rightarrow \T Q$ is the anchor map. In this case, the isotropy
$\mathfrak{h}_{Q}=\ker\rho$ is a locally trivial bundle of Lie algebras over $Q$ with typical fiber $\mathfrak{h}$. So, the co-isotropy $E_{Q}=\mathfrak{h}_{Q}^{\ast}$ becomes the Lie-Poisson bundle. There exists an exact sequence of vector bundles
    \begin{equation*}
        0\rightarrow\mathfrak{h}_{Q}\hookrightarrow A\overset{\rho}{\rightarrow}\T Q\rightarrow0.
    \end{equation*}
Every splitting $\gamma:\,\T Q\rightarrow A$ of this sequence induces a linear connection $\nabla$ on $E_{Q}$ satisfying conditions (\ref{LPVH1}%
)-(\ref{LPVH3}) \cite{Vo-04,Mac-95}.  In this case, $\mathbf{F}$ is given by
    \begin{equation*}
        \mathbf{F}(u_{1},u_{2}):=\{\gamma(u_{1}),\gamma(u_{2})\}_{A}-\gamma([u_{1},u_{2}])
    \end{equation*}
for $u_{1},u_{2}\in\mathfrak{X}(Q)$.
\end{example}

Remark that an important class of transitive Lie algebriods comes from principle bundles \cite{Mac-95}. As is known \cite{AlMo-85}, there are ``nonintegrable'' Lie algebroids, which are transitive and cannot be realized as the Lie algebroids of principle bundles.

Now suppose  we are given an action of a connected Lie group $G$ (not necessarily compact) on $E_{Q}$ with an infinitesimal generator of the form
    \begin{equation*}
        \Upsilon_{a}=\mathrm{ad}_{s_{a}}^{\ast}y\cdot\frac{\partial}{\partial y},
    \end{equation*}
where $\mathfrak{g}\ni a\mapsto s_{a}\in\Gamma(E_{Q}^{\ast})$ is a linear mapping. Then, this linear action is canonical relative to $\Psi$ with momentum map $\mathbb{J}$ given by $\mathbb{J}_{a}(y)=\langle s_{a}(\pi_{Q}(y)),y \rangle$.

We have the following criterion:

\begin{proposition}
For a pair $(\nabla,\mathbf{F})$ satisfying the conditions \ref{CC1}-\ref{CC3} above, the corresponding Poisson tensor $\Pi^{\nabla,\mathbf{F}}$ and the (generalized) Wong system \eqref{NW} are invariant with respect to the $G$-action if%
    \begin{equation}
        \mathrm{ad}_{\nabla^{\ast}s_{a}}=0 \quad \forall a\in\mathfrak{g}.\label{ICO}%
    \end{equation}
\end{proposition}

The proof of this proposition is based on the following identity
    \begin{equation*}
        [\mathrm{hor}^{\nabla}_{v},\Upsilon_{a}] = \mathrm{ad}^{\ast}_{\nabla^{\ast}s_{a}}y \cdot \frac{\partial}{\partial y}
    \end{equation*}
which says that condition \eqref{ICO} implies the $G$-invariance of the horizontal distribution $\mathbb{H}^{\nabla} \subset \T{E_{Q}}$.

\begin{corollary}
Suppose that $E_{Q} \to Q$ is a locally trivial Lie-Poisson bundle whose typical fiber is the coalgebra of a Lie algebra $\mathfrak{h}$. If $\mathfrak{h}$ is semi-simple, then the $G$-invariance condition \eqref{ICO} reads
    \begin{equation*}
        \nabla^{\ast}s_{a}=0 \quad \forall a \in \mathfrak{g}.
    \end{equation*}
\end{corollary}

Given sections $s_{a}$, one can think of condition \eqref{ICO} as an equation for a connection $\nabla$. Of course, under condition (\ref{ICO}), the Wong system (\ref{NW}) is also $G$-invariant.

In order to resolve (\ref{ICO}) for the action of a compact Lie group $G$,
%For the action of a compact Lie group $G$, to resolve ,
one can try the following strategy: fixing a pair $(\nabla_{0},\mathbf{F}_{0})$ satisfying conditions (\ref{LPVH1})-(\ref{LPVH3}), we apply the averaging procedure \cite{MaMoRa-90,JYu,MYu-2012,MHVV} to get some $G$-invariant data $(\nabla,\mathbf{F})$ preserving the desired properties. In this case, one can show that a \emph{$G$-invariant connection} $\nabla$ on $E_{Q}$ is defined as $\nabla := \nabla_{0} - \mathrm{ad}^{\ast}_{\mathbf{A}}$, where a gauge 1-form $\mathbf{A} \in \Omega^{1}(Q) \otimes \Gamma(E^{\ast}_{Q})$ is given by the following ``linear'' version of formula \eqref{IGF1}:
    \begin{equation}\label{eq:Alinear}
        \mathbf{A} = - \int_{G}\left( \int_{0}^{1} \mathrm{Ad}_{\exp(ts_{a})}\nabla_{0}^{\ast}s_{a} \dd{t}\right) \dd{g}, \quad g=\exp{a}.
    \end{equation}
Consider the corresponding 2-form $\mathbf{F} \in \Omega^{1}(Q) \otimes \Gamma(E_{Q}^{\ast})$ defined by
    \begin{equation}\label{eq:Flinear}
        \mathbf{F} = \mathbf{F}_{0} + \nabla_{0}\mathbf{A} + \frac{1}{2}[\mathbf{A} \wedge \mathbf{A}]_{\mathrm{fib}}.
    \end{equation}

\begin{proposition}
For a given pair $(\nabla_{0},\mathbf{F}_{0})$, formulas \eqref{eq:Alinear} and \eqref{eq:Flinear} give the $G$-invariant data $(\nabla,\mathbf{F})$ satisfying \eqref{LPVH1}-\eqref{LPVH3} and inducing the generalized Wong system \eqref{NW} with $G$-symmetry.
\end{proposition}

In the trivial case $E_{Q}=Q\times\mathfrak{h}^{\ast}$, for some given $s_{a}\in C^{\infty}(Q;\mathfrak{h)}$, equation (\ref{ICO}) for $\nabla=\nabla^{\mathbf{A}}$ takes the form
    \begin{equation}
        \mathrm{ad}_{s_{a}}\mathbf{A=-}\dd_{Q}s_{a} \quad \left(
        \mathrm{mod}\,\Omega^{1}(Q;\mathrm{Cent}(\mathfrak{h}))\right)
        \label{AE}%
    \end{equation}
for all $a\in\mathfrak{g}$. Starting with canonical connection $\nabla_{0}$ on $Q\times\mathfrak{h}^{\ast}$ with $\mathbf{F}_0 = 0$, a solution to \eqref{AE} is just given by \eqref{eq:Alinear}.

Finally, we conclude with the following alternative derivation of formula \eqref{eq:Aiqy}, by solving equation \eqref{AE} for $\mathfrak{h} = \mathrm{so}(3)$ and a given $s \in C^{\infty}(Q;\mathbb{R}^{3})$ with $|s(q)|= 1$. In this case, \eqref{AE} is written as
    \begin{equation*}
        s \times \mathbf{A} = \dd_{Q}s
    \end{equation*}
and has a solution $\mathbf{A}=-s \times \dd_{Q}s$.
%If the Lie group $G$ is compact, then solutions to these equations are given by formula (\ref{IGF1}). Here, according to the general strategy, we start with a trivial connection $\nabla_{0}$ on $Q\times\mathfrak{h}^{\ast}$ with $\mathbf{F}_{0}=0$.

%In the noncompact case, for example when the Lie algebra $\mathfrak{h}$ is semisimple, for solving (\ref{AE}) one can try to apply some algebraic arguments.

\paragraph{Acknowledgements.} This research was partially supported by the Mexican National Council of Humanities, Science and Technology (CONAHCyT) under the grant CF-2023-G-279 and the University of Sonora (UNISON) under the project no. USO315008535.


\begin{thebibliography}{99}

\bibitem{AlMo-85}  Almeida, R.;  Molino, P. \textit{Suites d'Atiyah et
feuilletages transversalement copmlets,} C. R. Acad. Sci. Paris Ser. I, Math.,
vol. 300 (1985), 13-15.

\bibitem {MHVV} Avenda\~{n}o-Camacho, M.; Hasse-Armengol, I.;  Velasco-Barreras, E.;
 Vorobjev, Yu. The method of averaging for Poisson connections on
foliations and its applications,\textit{ J. Geom. Mech.} \textbf{2020}, 12, 343-361.

\bibitem{MYu-2012} Avendano-Camacho, M.; Vorobiev, Yu. The averaging method on slow-fast spaces with symmetry. \textit{Journal of Physics: Conf. Ser.} \textbf{2012}, 343.

\bibitem {MYu-17} Avenda\~{n}o-Camacho, M.;   Vorobiev, Yu. Deformations of
Poisson structures on fibered manifolds and adiabatic slow-fast systems. \textit{ Int.
J. Geom. Methods Mod. Phys.} \textbf{2017}, 14, 15 pages.

\bibitem {GHV} Greub, W.;  Halperin, S.;  Vanstone, R. \textit{Connections,
curvature, and cohomology}, Vol.~II, Academic Press, New York--London, 1973.

\bibitem{Holten-07} van Holten, J. W. Covariant Hamiltonian dynamics. \textit{ Phys.
Rev. D } \textbf{2007}, 75, 025027.

\bibitem {HorvNG-06}  Horv\'athy, P. A.; Ngome, J. P. Conserved quantities in
non-abelian monopole fields. \textit{ Phys. Rev. D} \textbf{2009}, 79, 127701.

\bibitem {Land-98} Landsman, N. P. \textit{ Mathematical topics between classical
and quantum mechanics}. Springer-Verlag, N.Y., 1998.

\bibitem {Mac-95}  Mackenzie, K.C.H. \textit{ Lie algebroids and Lie pseudoalgebras},
Bull. London Math. Soc. 27, 1995; pp. 97-147.

\bibitem{MaMoRa-90} Marsden, J. E.;  Montgomery, R.;  Ratiu, T. \textit{Reduction,
symmetry, and phase in mechanics,} Memoirs of AMS, Providence, vol.88, no.436
1990; pp. 1-110.

\bibitem {Maspfuhl-04} Maspfuhl, O.
Wong's equations in Poisson Geometry. \textit{J. Symplectic Geom.}
\textbf{2004}, 2, 545-578

\bibitem{Mon-84} Montgomery, R. Canonical formulations of a classical particle
in a Yang-Mills Field and Wong's equations. \textit{Lett. in Math. Phys.} \textbf{1984}, 8,  59-67.

\bibitem {MoMaRa-84} Montgomery, R.; Marsden,  J. E.;  Ratiu, T. Gauged
Lie--Poisson structures. In: \textit{Fluids and Plasmas: Geometry and Dynamics}
(J. Marsden, ed.) Cont. Math., vol. 28, 1984; pp. 101--114.

\bibitem {OrtRa-04} Ortega, J. P.;  Ratiu, T. S. \textit{Momentum maps and
Hamiltonian reduction}. Progress in Math., v.222, Birkh\"{a}user Boston Inc.,
Boston, 2004.

\bibitem {Sternberg-1977} Sternberg, S.
On minimal coupling and symplectic mechanics of a classical particle in the presence of a Yang-Mills field.
\textit{Proc. Nat. Acad. Sci.} \textbf{1977}, 74, 5253–5254.

\bibitem {Va-98} Vaisman, I. \textit{Lectures on the geometry of Poisson
manifolds}, Progress in Math., vol.118, Birkh\"{a}user, Boston, 1994, 205 pp.

\bibitem{JYu} Vallejo, J. A.; and  Vorobev; Yu. Invariant Poisson realizations and
the averaging of Dirac structures. \textit{Symmetry, Integrability and Geometry:
Methods and Applications, SIGMA} \textbf{2014}, 096, 1-20

\bibitem {Vo-01} Vorobjev, Yu. Coupling tensors and Poisson geometry near a
single symplectic leaf. In: \textit{ Lie Algebroids and Related Topics in Differential
Geometry}, Banach Center Publ., vol. 54, Polish Acad. Sci., Warsaw, 2001, pp. 249-274

\bibitem {Vo-04} Vorobjev, Yu. On Poisson realizations of transitive Lie
algebroids.\textit{ J. of Nonlinear Math. Phys.} \textbf{2004}, 11, 43-48.

\bibitem {Weinstein-1978} Weinstein, A.
A universal phase space for particles in Yang-Mills fields.
\textit{Lett. Math. Phys.} \textbf{1978}, 2, 417–420.

\bibitem {Wong} Wong, S. K. Field and particle equations for the
classical Yang-Mills field and particles with isotopic spin.\textit{ Nuovo
Cimento} \textbf{1970}, 65,  689-694.

%%%%%%%%%%%%%%%%%%%%%%%%%%%%%%%%%%%%%%%%%%%%%%%%%%%%%%%%%%%%%%%%%%%%%%%%%%%%%%%%%%%%%%%%%%%%%%%%%%%%%%%%%%%%%%%%%%

%\bibitem {Ginz-1}V. L. Ginzburg. Momentum mappings and Poisson cohomology.
%\textit{Internat. J. Math.} \textbf{1996}, 7, 329--358.

\end{thebibliography}
\end{document}